\documentclass[letterpaper, 10 pt, conference]{ieeeconf}
\IEEEoverridecommandlockouts  % Needed only to use the \thanks command
\usepackage[utf8]{inputenc}
\usepackage[T1]{fontenc}
\usepackage[english]{babel} % creates problem with abstract???
\usepackage{cite}
\usepackage{graphicx}
\usepackage{times}
\usepackage{amsfonts,amsmath,amssymb} % assumes amsmath package installed
\usepackage{url}

\graphicspath{{Figures/}}

\newtheorem{proposition}{Proposition}
\newtheorem{remark}{Remark}

\DeclareMathOperator{\sat}{sat}
\newcommand{\abs}[1]{\left\lvert#1\right\rvert}

\newcommand{\Rset}{\mathbb{R}}

\title{\Large A semi-global model-based state observer for the quadrotor using only inertial measurements}
\author{Philippe~Martin and Ioannis~Sarras
\thanks{P. Martin and I. Sarras are with the Centre Automatique et Syst\`emes, Mines ParisTech, PSL Research University, Paris, France. (e-mail: philippe.martin@mines-paristech.fr, ioannis.sarras@mines-paristech.fr).}
\thanks{This work was supported by the French Agence Nationale de la
Recherche through the ANR ASTRID SCAR project ``Sensory Control of
Aerial Robots" (ANR-12-ASTR-0033).}}

\author{Philippe Martin and Ioannis Sarras% <-this % stops a space
\thanks{This work was supported by the French Agence Nationale de la Recherche through the ANR ASTRID SCAR project ``Sensory Control of
Aerial Robots'' (ANR-12-ASTR-0033)}% <-this % stops a space
\thanks{P.~Martin and I.~Sarras are with the Centre Automatique et Systèmes, MINES ParisTech, PSL Research University, 75006 Paris, France
{\tt\footnotesize \{philippe.martin,ioannis.sarras\}@mines-paristech.fr}}%
}

\begin{document}
\maketitle

\begin{abstract}
We propose a nonlinear observer to estimate the state (orientation and in-plane velocity vector) of the quadrotor, based on a drag-force-enhanced model. It is a simpler and more robust alternative to recent works using a similar model together with an Extended Kalman Filter (EKF). A particular state over-parameterization leads to a linear time-varying model with a nonlinear state-constraint that serves for the observer design. The proposed observer is able to ensure the uniform semi-global asymptotic stability of zero estimation error by incorporating the nonlinear constraint into the correction terms.
\end{abstract}

\section{Introduction}
When designing a stabilizing control law for a quadrotor, the estimation of the state (orientation and linear velocity vector) of the quadrotor is usually paramount. While many sensors can potentially be used, the ability to control the aircraft relying only on inertial sensors (namely strapdown MEMS accelerometers and rate gyros) remains an important issue. Indeed, the inertial sensors are not subjected to temporary outages and are not disturbed by environmental perturbations.

In control laws relying solely on inertial sensors, it is usually assumed that a strapdown triaxal accelerometer measures the gravity vector resolved in body coordinates. According to~\cite{Martin:2010}, this is a rather poor dynamic approximation; for that reason \cite{Martin:2010} introduces an enhanced model taking into account the so-called rotor drag. While the rotor drag terms are rather small, they are nevertheless first-order and dominant in the accelerometer measurements. As they give an image of the linear velocity in the quadrotor plane, they provide some very interesting information for control purposes. A simple control scheme based on a linear observer is also proposed~\cite{Martin:2010} to demonstrate the possible benefits of the approach. 

More elaborate observers relying on this enhanced model (with some simplifications) are proposed and thoroughly investigated in \cite{Abeywardena:2013} and~\cite{Leishman:2014}, using an Extended Kalman Filter (EKF). The EKF has the great advantage of a general and systematic method, but has also several drawbacks: the convergence of the estimation can usually be guaranteed only on slowly varying trajectories; it is not so easy to tune and to initialize; finally, it is computationally rather heavy for an embedded processor.

In this paper, which builds on the preliminary version \cite{MartiS2015arxiv}, we propose a deviation from the main literature that exploits (local) EKF-based estimators based on the enhanced model and we show that under the same assumptions as~\cite{Abeywardena:2013} and~\cite{Leishman:2014}, a simple nonlinear observer with constant gains can provide a semi-globally convergent state estimate. From a practical viewpoint, and in comparison to an EKF, the proposed observer ensures a guaranteed convergence even for aggressive trajectories, as well as a much lower computational cost and an easier gain tuning. Notice that a different approach but similar in spirit, i.e. a nonlinear observer based on an enhanced aerodynamic model, is proposed in~\cite{AllibABM2014CCA}.

An original aspect of the proposed observer is that it deliberately ``ignores'' the geometry of the system: the estimated orientation does not live on the unit sphere (it merely converges to the sphere), but instead in a linear space of a higher dimension. This idea of ``ignoring'' the geometry to benefit from a linear structure in a higher dimensional space is borrowed from the recent developments~\cite{Batista:2014a,Batista:2014b,Eudes:2014} where the estimation state space lives in a higher dimensional space and estimation is viewed as a linear (time-varying) problem. This amounts to the non-trivial problem of designing an observer for a time-varying linear system with a (nondifferential) nonlinear constraint. A second novelty of the proposed observer, whose design is based on a modification of the approach of~\cite{InIbook,Karagiannis:2009} and references therein, lies in the fact that the nonlinear constraint is directly incorporated into the correction term. This is the key ingredient for obtaining a strict Lyapunov function and establishing the uniform semi-global asymptotic stability claim.

The paper runs as follows: section~\ref{sec:Models} presents the drag-forced enhanced model and its simplified version; in section~\ref{sec:Observer} the observer is presented and its convergence is proved; section~\ref{sec:tuning} is devoted to the local tuning of the gains from the linearized error system; in section~\ref{sec:Simulations}, simulations illustrate the good behavior of the observer.

\section{Models}\label{sec:Models}
We consider a quadrotor equipped with strapdown triaxal rate gyro and accelerometer located at the center of mass. The rate gyro measures the angular velocity vector projected in body axes, while the accelerometer measures the specific acceleration projected in body axes. To design a state estimator, we need a model of the dynamics and of the quantities measured by the sensors.

\subsection{The drag-force-enhanced model of~\cite{Martin:2010}}
In the drag-force-enhanced model introduced in~\cite{Martin:2010}, the translation dynamics and orientation kinematics of the quadrotor are described by
\begin{IEEEeqnarray}{rCl}
	\dot u &=& vr-wq+g\eta_1-cu\label{eq:ufull}\\
	\dot v &=& wp-ur+g\eta_2-cv\label{eq:vfull}\\
	\dot w &=& uq-vp+g\eta_3-\frac{T}{m}\\
	\dot\eta_1 &=& r\eta_2-q\eta_3\label{eq:eta1full}\\
	\dot\eta_2 &=& p\eta_3-r\eta_1\\
	\dot\eta_3 &=& q\eta_1-p\eta_2\label{eq:eta3full}\\
	1 &=& \eta_1^2+\eta_2^2+\eta_3^2,\label{eq:constraintfull}
\end{IEEEeqnarray}
where $u$, $v$ and $w$ are the coordinates in body axes of the mass center velocity vector; $p$, $q$, $r$ are the coordinates in body axes of the angular velocity vector; $\eta_1$, $\eta_2$ and $\eta_3$ are the entries of the third column of the rotation matrix from the Earth frame to the body frame; $T$ is the total thrust, $g$ is the gravity constant, and $m$ is the quadrotor mass. Finally, the (time-varying) drag coefficient~$c$ reads
\begin{IEEEeqnarray*}{rCl}
	c(t) &:=& \frac{\lambda}{m}\sum_{i=1}^{4}\omega_i(t),
\end{IEEEeqnarray*}
where $\lambda>0$ is constant, and $\omega_i\ge0$ is the $i$-th motor speed.

The constraint~\eqref{eq:constraintfull} is just the fact that the sum of the squares of every column (and line) of a rotation matrix is~$1$. When the rotation matrix is parametrized by the roll, pitch and yaw angles $(\phi,\theta,\psi)$, we have
\begin{IEEEeqnarray*}{rCl}
	(\eta_1,\eta_2,\eta_3) &=& (-\sin\theta,\sin\phi\cos\theta,\cos\phi\cos\theta),
\end{IEEEeqnarray*}
and the constraint~\eqref{eq:constraintfull} is automatically enforced. The kinematic equations~\eqref{eq:eta1full}--\eqref{eq:eta3full} then amount to
\begin{IEEEeqnarray*}{rCl}
	\dot\phi &=&  p + (q\sin\phi + r\cos\phi)\tan\theta\\
	\dot\theta &=& q\cos\phi-r\sin\phi.
\end{IEEEeqnarray*}

The difference between the enhanced and the standard model lies in the drag terms $cu,cv$ in \eqref{eq:ufull}-\eqref{eq:vfull}, due to the so-called rotor drag. While these terms are rather small, they are nevertheless dominant in the accelerometer measurements. Indeed the triaxial accelerometer measures
\begin{IEEEeqnarray*}{rCl}
	(a_x,a_y,a_z) &=&  \Bigl(-cu,-cv,-\frac{T}{m}\Bigr),
\end{IEEEeqnarray*}
which is the specific acceleration vector (i.e. all the external forces except gravity). This model of the accelerometer measurements is much more accurate and interesting for control purposes than the often used rather crude approximation
\begin{IEEEeqnarray*}{rCl}
	(a_x,a_y,a_z) &=&  (g\sin\theta,-g\sin\phi\cos\theta,-g\cos\phi\cos\theta).
\end{IEEEeqnarray*}
Finally, the rate gyro provides the measurements
\begin{IEEEeqnarray*}{rCl}
	(g_x,g_y,g_z) &=&  (p,q,r).
\end{IEEEeqnarray*}

\subsection{A simplified design model}
To be able to design an estimator relying only on inertial measurements, we make the following assumptions:
\begin{itemize}
	\item the Coriolis terms are small enough to be neglected
	\item the time-varying coefficient $c$ is known, with known constant bounds $0\le c_l\le c(t) \le c_u$
	\item the measurements from the accelerometer and the rate gyro are perfect, i.e. unbiased and noiseless.
	%\item no wind perturbation is considered.
\end{itemize}
The first assumption, also used in~\cite{Abeywardena:2013,Leishman:2014}, is reasonable except for very aggressive trajectories. The second assumption is satisfied since the motor speeds are usually available and the constant $\lambda$ can be experimentally evaluated; alternatively $c$ can be replaced by its nominal value in hovering $\bar{c}:=\frac{4\lambda\bar{\omega}}{m}$, as in~\cite{Abeywardena:2013,Leishman:2014}. The robustness of the observer with respect to the neglected Coriolis forces and unaccounted biases and noises is illustrated in simulation in section \ref{sec:Simulations}.

Thanks to the first assumption, the vertical velocity~$w$ does not influence the other equations, so that we do not need to consider its evolution. The simplified model we will use to design the observer is therefore
\begin{IEEEeqnarray}{rCl}
	\dot u &=& g\eta_1-cu\label{eq:u}\\
	\dot v &=& g\eta_2-cv\\
	\dot\eta_1 &=& r\eta_2-q\eta_3\\
	\dot\eta_2 &=& p\eta_3-r\eta_1\\
	\dot\eta_3 &=& q\eta_1-p\eta_2\\
	1 &=& \eta_1^2+\eta_2^2+\eta_3^2,\label{eq:constraint}
\end{IEEEeqnarray}
where $u,v,p,q,r$ are known since measured (assuming perfect sensors). This is the same design model as in~\cite{Abeywardena:2013,Leishman:2014}, except that we stick to the~$\eta_1,\eta_2,\eta_3$ variables instead of using the parametrization by the roll and pitch angles, and that no gyro bias is considered. 
A benefit of the representation used here is that it yields a (time-varying) linear differential system with a constraint; we will take advantage of this structure by relaxing the constraint in the state estimator, and using it rather as a feedback signal. The linear structure then makes the design of the estimator and the proof of convergence relatively easy.

Finally, notice the system~\eqref{eq:u}--\eqref{eq:constraint} is observable, meaning all the state variables can be expressed as functions of the measured quantities and their derivatives, provided $\eta_3$ remains strictly positive. Indeed, we then have
\begin{IEEEeqnarray*}{rCl}
	\eta_1 &=&  \frac{\dot u+cu}{g}\\
	\eta_2 &=&  \frac{\dot v+cv}{g}\\
	\eta_3 &=& \sqrt{1-\eta_1^2-\eta_2^2}.
\end{IEEEeqnarray*}
The condition $\eta_3>0$ is satisfied if the roll and pitch angles remain within~$(-\pi/2,\pi/2)$, which is the case except for extremely aggressive trajectories.

\section{A semi-global observer}\label{sec:Observer}
In this section, we show we can estimate the state of the design system~\eqref{eq:u}--\eqref{eq:constraint} by the system
\begin{IEEEeqnarray}{rCl}
	\dot{\hat{u}} &=& g\hat{\eta}_1-c\hat{u} - (k_u+k_1)(\hat{u} + \frac{a_x}{c})\label{eq:uhat}\\
	\dot{\hat{v}} &=& g\hat{\eta}_2-c\hat{v} -(k_v+k_2)(\hat v + \frac{a_y}{c})\\
	\dot{\hat\eta}_1 &=& r\hat{\eta}_2 - q\hat{\eta}_3 - \frac{k_1 k_u}{g}(\hat{u} + \frac{a_x}{c}) -\frac{rk_2}{g}(\hat v + \frac{a_y}{c})\\
	\dot{\hat\eta}_2 &=&  p\hat{\eta}_3-r\hat{\eta}_1 +\frac{rk_1}{g}(\hat u + \frac{a_x}{c}) - \frac{k_2 k_v}{g}(\hat{v} + \frac{a_y}{c})\\
	\dot{\hat\eta}_3 &=& q\hat{\eta}_1 - p\hat{\eta}_2 -\frac{qk_1}{g}(\hat u+\frac{a_x}{c})+\frac{pk_2}{g}(\hat v+\frac{a_y}{c})\\
	&& \,-k_3E\Bigl(\hat{\eta}_1-\frac{k_1}{g}\bigl(\hat{u}+\frac{a_x}{c}\bigr),\hat{\eta}_2-\frac{k_2}{g}\bigl(\hat v+\frac{a_y}{c}\bigr),\hat{\eta}_3\Bigr),\IEEEeqnarraynumspace\label{eq:eta3hat}
\end{IEEEeqnarray}
where
\begin{IEEEeqnarray*}{rCl}
	E(x_1,x_2,x_3) &:=& \frac{\displaystyle x_3 - \sqrt{1 - \sat(x_1^2 + x_2^2)}}{\displaystyle x_1^2+x_2^2 + 1 - \varepsilon^2}\\
	\sat(x) &:=& \min(1,\frac{1-\varepsilon^2}{|x|})x;
\end{IEEEeqnarray*} 
$\varepsilon>0$ is a sufficiently small constant, while $k_1,k_2,k_3,k_u,k_v$ are (constant) tuning gains yet to be chosen. Notice~\eqref{eq:uhat}--\eqref{eq:eta3hat} has the classical structure of an observer, namely a copy of the design system plus correction terms which are zero when the estimated state equals the actual state.

\begin{proposition}
The system ~\eqref{eq:uhat}--\eqref{eq:eta3hat} is a (uniformly) semi-globally asymptotically convergent observer of the design system~\eqref{eq:u}--\eqref{eq:constraint}. In other words: assume the considered trajectory of~\eqref{eq:u}--\eqref{eq:constraint} satisfies $\eta_3(t)\ge\varepsilon$ for all $t\ge0$; then 
\begin{IEEEeqnarray*}{l}
	\lim_{t\rightarrow+\infty}\bigl(\hat u(t),\hat v(t),\hat\eta_1(t),\hat\eta_2(t),\hat\eta_3(t)\bigr)=\\
	\hspace{3cm}\bigl(u(t),v(t),\eta_1(t),\eta_2(t),\eta_3(t)\bigr)
\end{IEEEeqnarray*} 
whatever the initial condition $\bigl(\hat u(0),\hat v(0),\hat\eta_1(0),\hat\eta_2(0),$ $\hat\eta_3(0)\bigr)$, provided the gains satisfy
\begin{IEEEeqnarray}{rCl}
k_3 &>& 0 \label{eq:k3}\\
k_1 &>& 1+\frac{k_3}{2\varepsilon^2}\\
k_2 &>& 1+\frac{k_3}{2\varepsilon^2}\\
k_u &>& \frac{k_1^2 c_u^2}{2g^2}+\frac{g^2}{2}\\
k_v &>& \frac{k_2^2 c_u^2}{2g^2}+\frac{g^2}{2}\label{eq:kv}.
\end{IEEEeqnarray} 
\end{proposition}

\begin{remark}
	The qualifier ``semi-global'' refers to the fact that the domain of attraction and the gains depend on the choice of the parameter~$\varepsilon$. The smaller $\varepsilon$, the larger the domain of attraction, but the larger the gains $k_u,k_v,k_1,k_2$.
\end{remark}

\begin{proof}
Define the error variables
\begin{IEEEeqnarray*}{rCl}
e_u &:=& \hat u - u\\
e_v &:=& \hat v - v\\
z_1 &:=&  \hat\eta_1 - \eta_1 - \frac{k_1}{g}(\hat u - u)\\
z_2 &:=&  \hat\eta_2 - \eta_2 - \frac{k_2}{g}(\hat v - v)\\
z_3 &:=&  \hat\eta_3 - \eta_3. 
\end{IEEEeqnarray*}
The error system reads
\begin{IEEEeqnarray*}{rCl}
\dot e_u &=& -(c+k_u)e_u+gz_1\label{eq:eu}\\
\dot e_v &=& -(c+k_v)e_v+gz_2\label{eq:ev}\\
\dot z_1 &=& rz_2-qz_3-k_1z_1+\frac{k_1c}{g}e_u\\
\dot z_2 &=& pz_3-rz_1-k_2z_2+\frac{k_2c}{g}e_v\\
\dot z_3 &=& qz_1-pz_2-k_3E(\eta_1+z_1,\eta_2+z_2,\eta_3+z_3).
\end{IEEEeqnarray*}

Consider now the quadratic function
\begin{IEEEeqnarray*}{rCl}
V &:=& \frac{1}{2}(e_u^2 + e_v^2)+Z, \label{eq:V}
\end{IEEEeqnarray*}
where $Z:=\frac{1}{2}(z_1^2+z_2^2+z_3^2)$. Using $ab\le\frac{a^2}{2}+\frac{b^2}{2}$ yields
\begin{IEEEeqnarray*}{rCl}
\frac{1}{2}\frac{d}{dt}e_u^2 &=& ge_uz_1-(c+k_u)e_u^2\\
&\le& \frac{z_1^2}{2}-\Bigl(c_l+k_u-\frac{g^2}{2}\Bigr)e_u^2\\
\frac{1}{2}\frac{d}{dt}e_v^2 &\le& \frac{z_2^2}{2}-\Bigl(c_l+k_v-\frac{g^2}{2}\Bigr)e_v^2,
\end{IEEEeqnarray*}
In the same way,
\begin{IEEEeqnarray*}{rCl}
\dot Z &=& \frac{k_1c}{g}e_uz_1-k_1z_1^2+\frac{k_2c}{g}e_vz_2-k_2z_2^2-k_3z_3E\\
&\le& \frac{z_1^2}{2}-k_1z_1^2+\frac{z_2^2}{2}-k_2z_2^2+\frac{k_1^2c_u^2}{2g^2}e_u^2+\frac{k_2^2c_u^2}{2g^2}e_v^2-k_3z_3E,
\end{IEEEeqnarray*}
where for brevity $E$ stands for $E(\eta_1+z_1,\eta_2+z_2,\eta_3+z_3)$.

We next bound the term $-k_3z_3E$. First, notice the function $F(x):=\sqrt{1-\sat(x)}$ satisfies
\begin{IEEEeqnarray*}{rCl}
\abs{F(x)-F(y)} &\le& \frac{1}{2\varepsilon}\abs{\sat(x)-\sat(y)}\\
&\le&  \frac{1}{2\varepsilon}\abs{x-y}.
\end{IEEEeqnarray*}
Set $A:=\eta_3-F\bigl((\eta_1+z_1)^2+(\eta_2+z_2)^2\bigr)$. The previous inequality and the assumption $\eta_3\ge\varepsilon$ yield
\begin{IEEEeqnarray*}{rCl}
\abs{A} &=& \abs{F\bigl(\eta_1^2+\eta_2^2\bigr)-F\bigl((\eta_1+z_1)^2+(\eta_2+z_2)^2\bigr)}\\
 &\le&  \frac{1}{2\varepsilon}\abs{(\eta_1+z_1)^2+(\eta_2+z_2)^2-\eta_1^2-\eta_2^2}\\
 &=& \frac{1}{2\varepsilon}\abs{z_1(2\eta_1+z_1)+z_2(2\eta_2+z_2)}.
\end{IEEEeqnarray*}
Using repeatedly $(a+b)^2\le2(a^2+b^2)$, this implies
\begin{IEEEeqnarray*}{rCl}
A^2 &\le&  \frac{2}{4\varepsilon^2}\bigl(z_1^2(2\eta_1+z_1)^2+z_2^2(2\eta_2+z_2)^2\bigr)\\
&\le& \frac{z_1^2}{\varepsilon^2}\bigl((\eta_1+z_1)^2+\eta_1^2\bigr) + \frac{z_2^2}{\varepsilon^2}\bigl((\eta_2+z_2)^2+\eta_2^2\bigr)\\
&\le& \frac{z_1^2}{\varepsilon^2}\bigl((\eta_1+z_1)^2+1-\varepsilon^2\bigr) + \frac{z_2^2}{\varepsilon^2}\bigl((\eta_2+z_2)^2+1-\varepsilon^2\bigr).
\end{IEEEeqnarray*}
In the last inequality, we have used 
\begin{IEEEeqnarray*}{C}
\eta_1^2,\eta_2^2\leq\eta_1^2+\eta_2^2=1-\eta_3^2\leq1-\varepsilon^2.
\end{IEEEeqnarray*}
Finally, setting $B:=(\eta_1+z_1)^2+(\eta_2+z_2)^2+1-\varepsilon^2$, and noticing
\begin{IEEEeqnarray*}{rCl}
B &\leq& 2(z_1^2+z_2^2) + 2(\eta_1^2+\eta_2^2)+ 1 - \epsilon^2\\
	&\leq& 2(z_1^2+z_2^2) + 3(1-\varepsilon^2)\\
	&\leq& 4V + 3(1-\varepsilon^2),
\end{IEEEeqnarray*}
we find
\begin{IEEEeqnarray*}{rCl}
-z_3E &=& -z_3\frac{z_3+A}{B}\\
&\le& -\frac{z_3^2}{B}+\frac{z_3^2+A^2}{2B}\\
&\le& -\frac{z_3^2}{2B}+\frac{z_1^2}{2\varepsilon^2}+\frac{z_2^2}{2\varepsilon^2}\\
&\le& -\frac{z_3^2}{8V+6(1-\varepsilon^2)}+\frac{z_1^2}{2\varepsilon^2}+\frac{z_2^2}{2\varepsilon^2}.
\end{IEEEeqnarray*}

We then collect all the previous findings to get
\begin{IEEEeqnarray*}{rCl}
\dot V &\le& -k_{10}z_1^2 - k_{20} z_2^2 - \frac{k_3z_3^2}{8V+6(1-\varepsilon^2)}- k_{u0}e_u^2 - k_{v0}e_v^2,
\end{IEEEeqnarray*}
where we have set
\begin{IEEEeqnarray*}{rCl}
k_u &:=&  k_{u0} + \frac{k_1^2 c_u^2}{2g^2}+\frac{g^2}{2}\\
k_v &:=&  k_{v0} + \frac{k_1^2 c_u^2}{2g^2}+\frac{g^2}{2}\\
k_1 &:=&  k_{10}+1+\frac{k_3}{2\varepsilon^2}\\
k_2 &:=&  k_{20}+1+\frac{k_3}{2\varepsilon^2}.
\end{IEEEeqnarray*} 
Notice $k_{u0},k_{v0},k_{10},k_{20}$ are by assumption strictly positive.
We conclude the proof by noticing that the function
\begin{IEEEeqnarray*}{rCl}
W &:=& 3(1-\varepsilon^2)V + 2V^2
\end{IEEEeqnarray*}
satisfies
\begin{IEEEeqnarray*}{rCl}
\dot W &\leq& -k_{10}z_1^2 - k_{20} z_2^2 - \frac{k_{3}}{2}z_3^2 -k_{u0}e_u^2 - k_{v0}e_v^2.
\end{IEEEeqnarray*}
We have assumed $3(1-\epsilon^2)\ge1$, which is clearly not restrictive.
%
%This concludes the proof of \emph{uniform (semi-global) asymptotic stability} of $(z,e_u,e_v)=(0,0,0)$ by incorporating Theorem 4.9 of \cite{Khalil} or Lemma 2.1 of \cite{MazencBook}.
\end{proof}

\begin{remark}
The saturation function defined above appears also in the context of quadrotor attitude estimation in \cite{Hua:2014a}. Furthermore, the parametrization $(\eta_1,\eta_2,\eta_3)$ has been exploited also in \cite{Hua:2014b} but with an observer respecting the inherent geometry of the system.
\end{remark}

\begin{remark}
The proposed observer was derived using the observer methodology based on invariant manifolds \cite{InIbook,Karagiannis:2009}. The general principle behind  this technique is to estimate the unmeasured state $\eta$ by rendering a certain manifold 
$$
\mathcal{M}=\{(\eta,y, \xi) \rvert \beta(\xi,y)=\varphi(\eta,y) \}
$$
attractive and invariant for some $\xi$ (the observer state), $y$ the measured state, and functions $\beta$, $\varphi$. The objective then is to stabilize to zero the dynamics of the ``error'' (usually called off-the-manifold coordinates)
$$z:=\beta(\xi,y) - \varphi(\eta,y),$$
whose norm essentially captures the distance from the manifold $\mathcal{M}$.
If this (non-standard) stabilization objective is achieved then an estimate of $\eta$ is given by $\hat{\eta} = \varphi^{-1}(\beta(\xi,y),y)$.

For our problem, we initially defined the mappings
\begin{IEEEeqnarray*}{rCl}
\beta(\xi,y)&:=& \xi + \beta_0(y)\\
\varphi(\eta,y) &:=& \eta
\end{IEEEeqnarray*}
$\beta_0(y)={\rm col}(\frac{k_1u}{g},\frac{k_2v}{g},0)$ that gave the error variable $z := \xi + \beta_0(y) - \eta$ with $\xi\in\Rset^3$ the state of the observer (to be designed) and $\beta_0:\Rset^2\to\Rset^3$ a mapping that was chosen properly. Our objective was then to show that $z$ is asymptotically stable with respect to the origin which implied that an estimate of $\eta$ was provided by $\hat{\eta} = \xi + \beta_0(y)$. However, since in practice it is preferable to use filtered versions of the measurements $\hat y=(\hat u, \hat v)$, we thus re-defined $\beta$ as $\beta(\xi,y, \hat{y}):= \xi + \beta_0(\hat y)= \xi + \beta_0(y) - (\beta_0(y)-\beta_0(\hat y))$ while having to additionally ensure that $\hat y - y$ was converging to zero.
The desired objective was then attained by incorporating a key additional term in the observer dynamics that involves the nonlinear state-constraint, and instead of following the usual design that hinges only upon the form of the system dynamics and the freedom on the mapping $\beta_0$, and which is crucial for obtaining a strict Lyapunov function.
\end{remark}

\section{Linearized system and gain tuning}\label{sec:tuning}
In the previous section, the semi-global convergence has been presented. However, it is also important in practice to ensure a good local behavior through a proper gain tuning. This is done by examining the linearized error system.

The first-order approximation of the error system in hovering, i.e. $(\bar{z},\bar{e}_u,\bar{e}_v):=(0,0,0)$, $\eta=(0,0,1)^T$, $(p,q,r):=(0,0,0)$, can be decomposed into the three decoupled subsystems
\begin{IEEEeqnarray}{rcl}
\left[ \begin{array}{c}
\dot{\zeta}_1\\
\dot{e}_u\\
\end{array}\right]
&=&\left[ \begin{array}{cc}
-k_1 & \frac{{\bar{c}k_1}}{g}\\
g & -(k_u+\bar{c})
\end{array}\right]
\left[ \begin{array}{c}
\zeta_1\\
e_u
\end{array}\right]\label{subsys1}
\end{IEEEeqnarray} 

\begin{IEEEeqnarray}{rCl}
\left[ \begin{array}{c}
\dot{\zeta}_2\\
\dot{e}_v\\
\end{array}\right]
&=&\left[ \begin{array}{cc}
-k_2 & \frac{{\bar{c}k_2}}{g}\\
g & -(k_v+\bar{c})
\end{array}\right]
\left[ \begin{array}{c}
\zeta_2\\
e_v
\end{array}\right]\label{subsys2}
\end{IEEEeqnarray} 

\begin{IEEEeqnarray}{rCl}
\dot{\zeta_3} &=& - \frac{k_{3}}{1-\epsilon^2}\zeta_3 \label{subsys3},
\end{IEEEeqnarray}
where $\bar{c}$ is the nominal value of the drag coefficient.

The characteristic polynomial of the first subsystem is $\lambda^2+(k_1+k_u+c)\lambda+k_1k_u\approx (\lambda+k_1)(\lambda+k_u)$, using the fact that $c\ll k_1+k_u$. We can hence deduce that the corresponding eigenvalues are approximately $\lambda=-k_1$, $\lambda=-k_u$. Similarly, the characteristic polynomial of the second subsystem is $\lambda^2+(k_2+k_v+c)\lambda+k_2k_v\approx (\lambda+k_2)(\lambda+k_v)$. The eigenvalue of the third subsystem is obviously $\lambda=-\frac{k_{3}}{1-\epsilon^2}$.

In conclusion, the linearized system has real negative eigenvalues that can be freely assigned through the choice of  $k_1,k_2,k_3,k_u,k_v>0$, provided the conditions~\eqref{eq:k3}-\eqref{eq:kv} are satisfied.

\section{Simulations}\label{sec:Simulations}

We illustrate the good behavior of the observer in simulation. The quadrotor is made to follow a rather aggressive trajectory. The robustness with respect to the design assumptions of neglected Coriolis forces and perfect measurements is investigated.

To this end, we consider that the measurements provided by the accelerometer and the rate gyro are corrupted by noise and constant biases, i.e.
\begin{IEEEeqnarray*}{rCl}
	(a_x,a_y) &=&  (-cu+b_u+\nu_u,-cv+b_v+\nu_v)\\
	(g_x,g_y,g_z) &=&  (p+b_p+\nu_p,q+b_q+\nu_q,r+b_r+\nu_r), 
\end{IEEEeqnarray*}
with $b_u,b_v,b_p,b_q,b_r$ and  $\nu_u,\nu_v,\nu_p,\nu_q,\nu_r$ the corresponding biases and noises.
The values of the biases are 
\begin{IEEEeqnarray*}{rCl}
	(b_u, b_v) &=&  (0.05,0.04)\quad(m/s^2),\\
	(b_p, b_q, b_r) &=& (0.02,-0.015,0.01)\quad (rad/sec),
\end{IEEEeqnarray*}
the noises are independent band-limited gaussian white noises with a noise power of $10^{-5}$ and sample time $10^{-4}$, filtered at $1kHz$.

The observer gains are set to
\begin{IEEEeqnarray*}{rCl}
	(k_1, k_2, k_3, k_u,k_v) &=&  (7,7,0.1,49,49),
\end{IEEEeqnarray*}
which corresponds to eigenvalues $-7$ and $-49$ for the two subsystems~\eqref{subsys1},~\eqref{subsys2}, and $-0.1$ for the subsystem~\eqref{subsys3}, assuming a nominal value $\bar{c}=0.25$. The very slow eigenvalue $-0.1$ is so chosen to avoid amplifying measurement noises, larger values can be selected for faster convergence.

\begin{figure}[h]
	\includegraphics[width=\columnwidth]{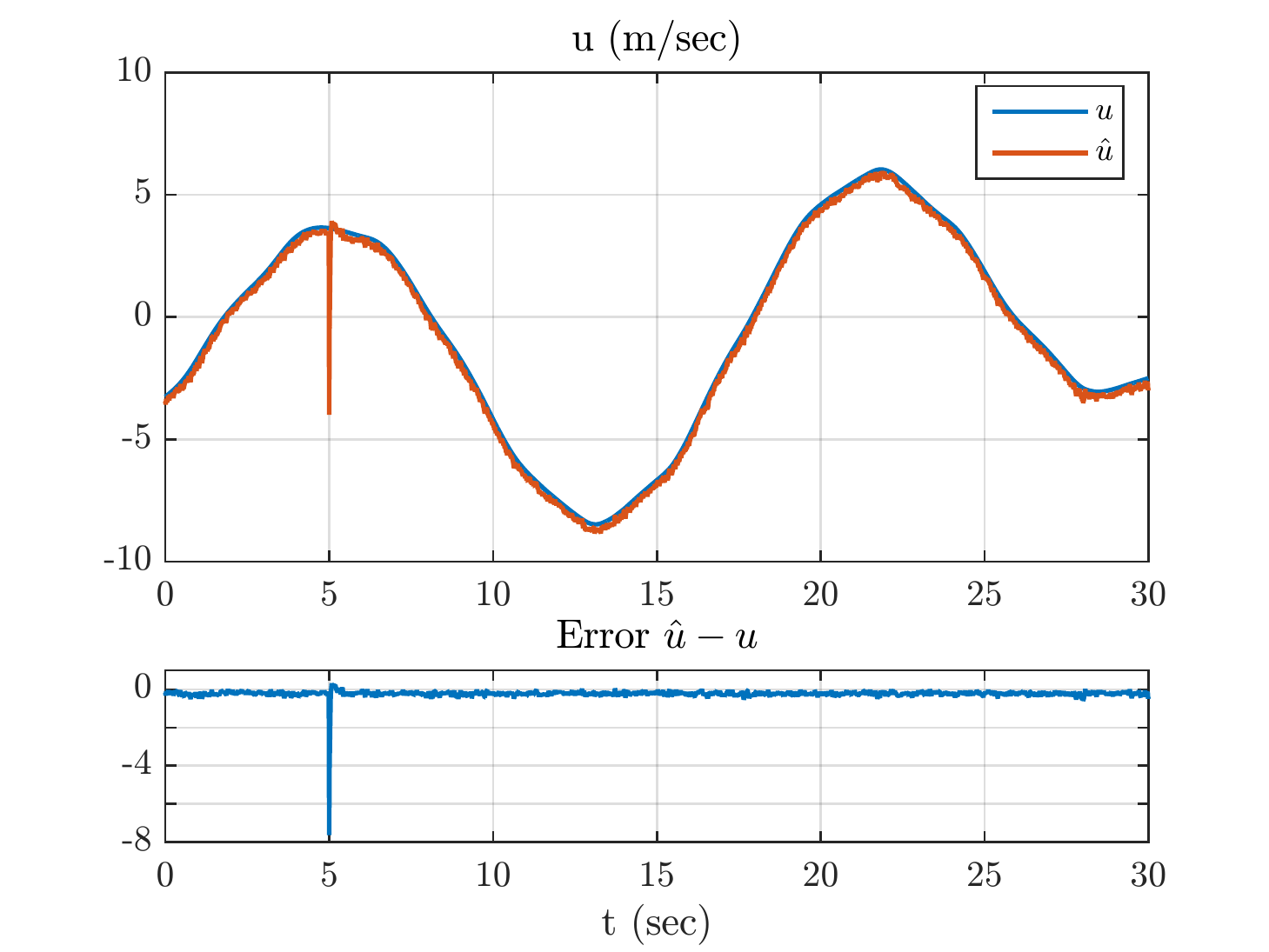}
	\caption{True (blue), estimated (red) velocity $u$ and error $\hat u -u$ ($m/sec$)}
	\label{fig:u}
\end{figure}

\begin{figure}[h]
	\includegraphics[width=\columnwidth]{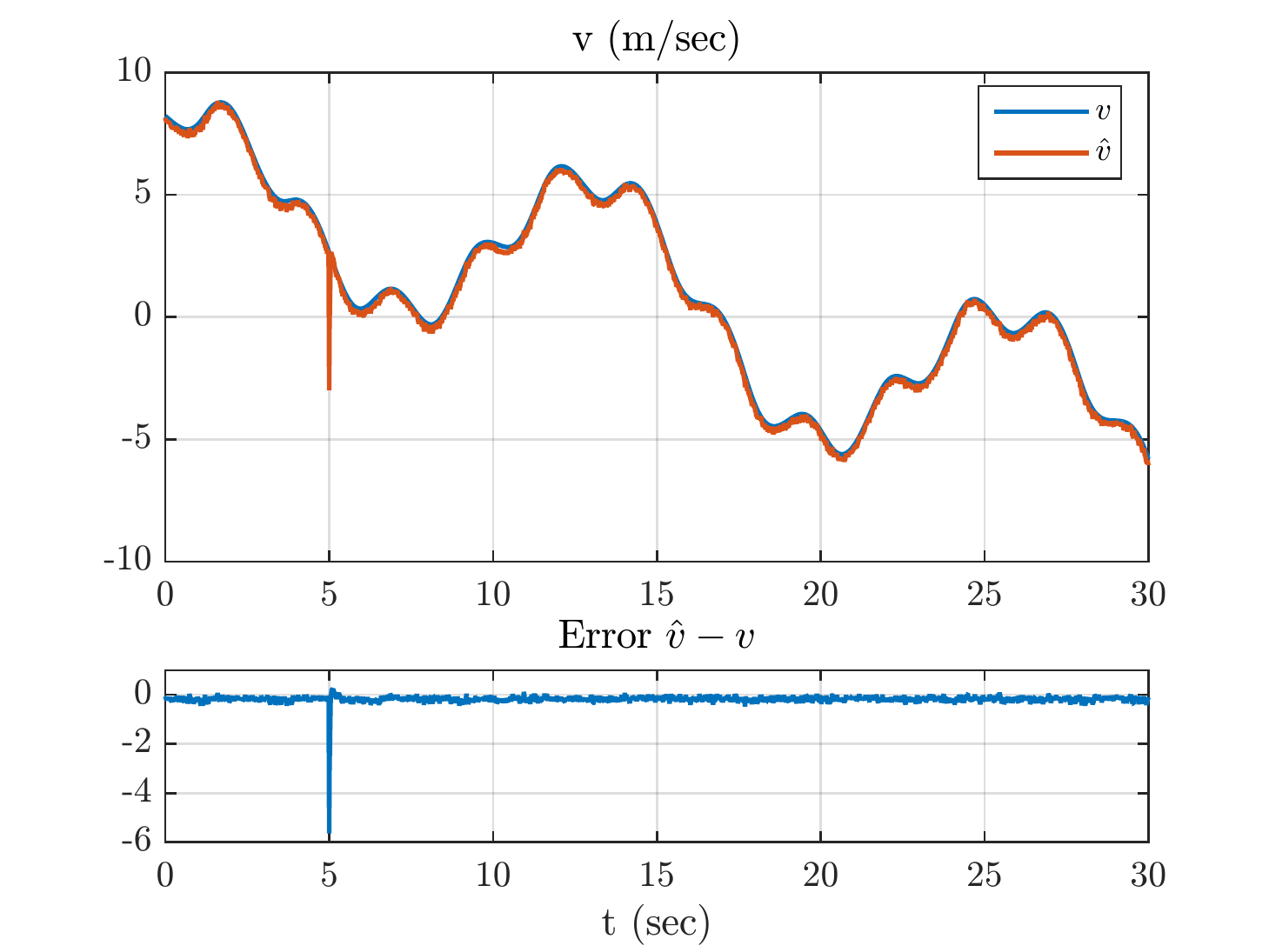}
	\caption{True (blue), estimated (red) velocity $v$ and error $\hat v- v$ ($m/sec$)}
	\label{fig:v}
\end{figure}

\begin{figure}[h]
	\includegraphics[width=\columnwidth]{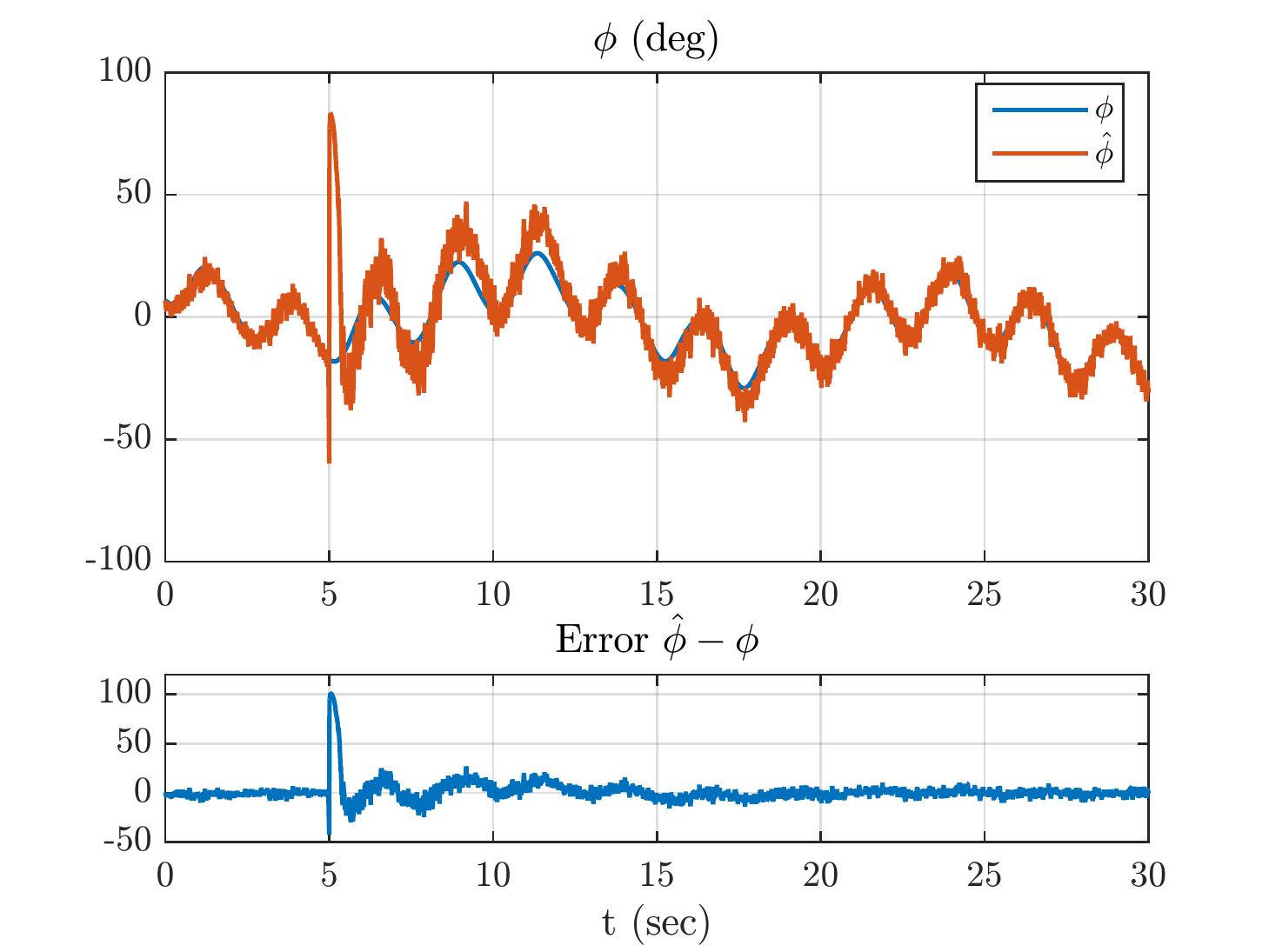}
	\caption{True (blue), estimated (red) $\phi$ and error $\hat\phi - \phi$ ($\deg$).}
	\label{fig:phi}
\end{figure}

\begin{figure}[h]
	\includegraphics[width=\columnwidth]{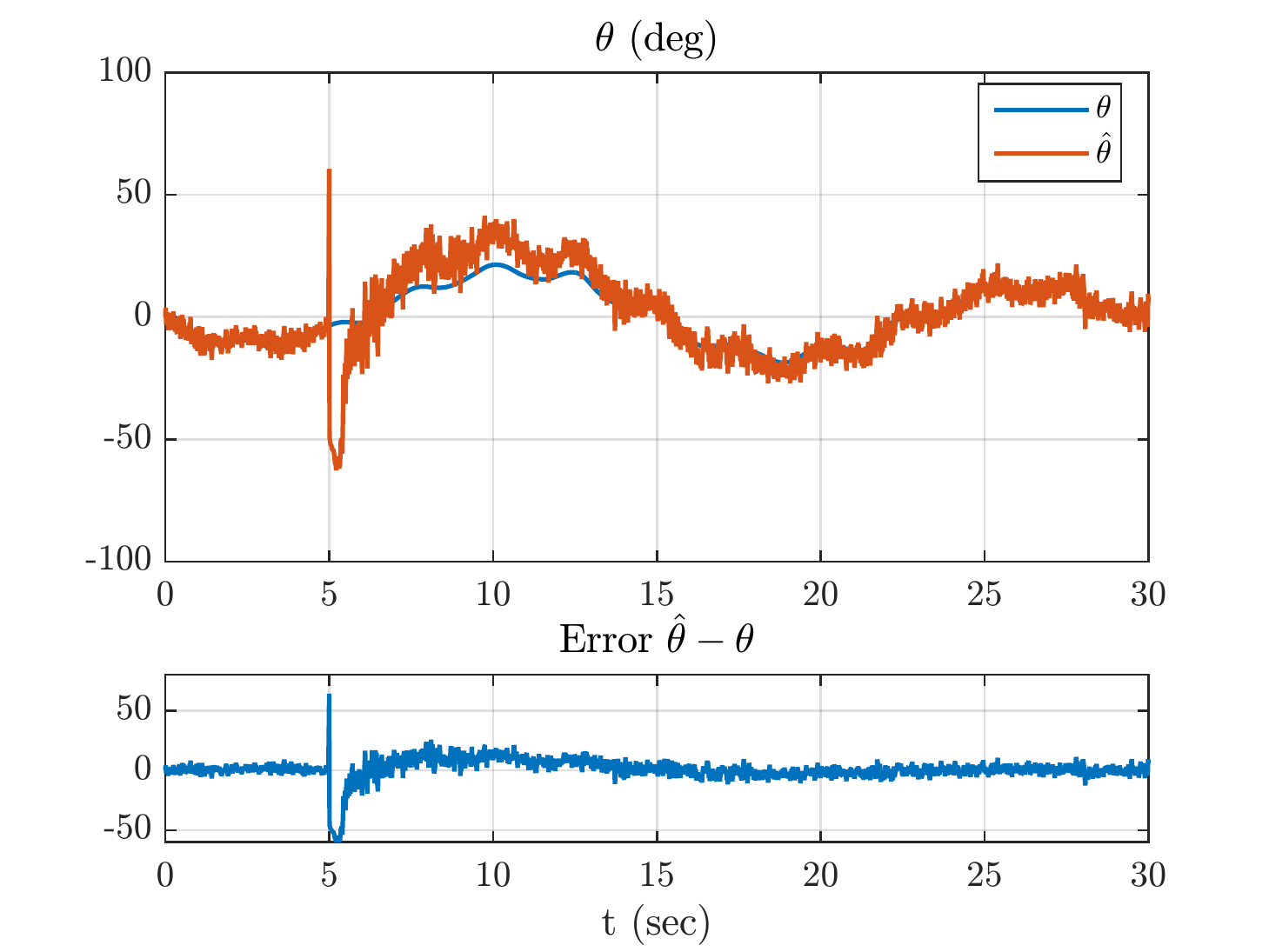}
	\caption{True (blue), estimated (red) $\theta$ and error $\hat\theta - \theta$ ($\deg$)}
	\label{fig:theta}
\end{figure}

\begin{figure}[h]
	\includegraphics[width=\columnwidth]{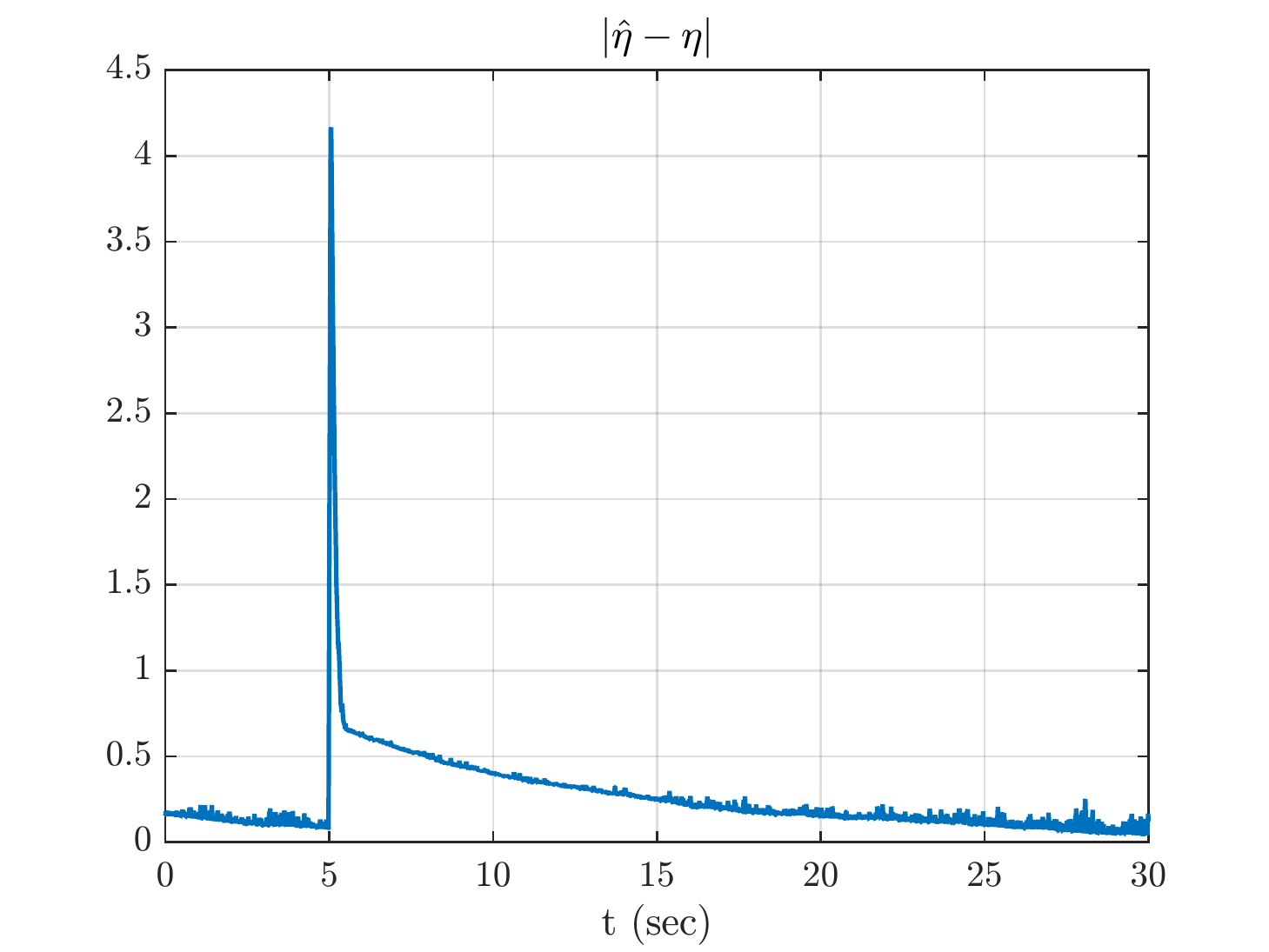}
	\caption{Norm of the attitude error $\hat\eta - \eta$.}
	\label{fig:error_eta}
\end{figure}

The observer is initialized with approximately no error, but is suddenly reinitialized  at the time instant $t=5~sec$ to $(\hat u(5) , \hat v(5), \hat \phi(5), \hat \theta(5)) = (-4~m/sec,-3~m/sec, -60~deg,60~deg)$.
The convergence of the estimated variables $\hat u,\hat v,\hat{\phi},\hat{\theta}$ is as anticipated excellent after the reinitialization, very fast for $\hat u,\hat v$ and slower for $\hat{\phi},\hat{\theta}$ (because of the slow eigenvalue $-0.1$), in accordance with the choice of gains, see Fig.\ref{fig:u}-\ref{fig:theta}. The estimated roll and pitch angles are obtained from $\hat{\eta}$ by $\hat\phi:=-\arcsin(\frac{\hat{\eta_1}}{|\hat{\eta}|})$ and $\hat\theta:=\arctan(\frac{\hat{\eta_2}}{\hat{\eta_3}})$. The norm of the attitude estimation error is also displayed in Fig.~\ref{fig:error_eta};  as $\hat\eta$ lives in $\mathbb{R}^3$, and not in $\mathtt{S}^2$ as $\eta$, hence $|\hat{\eta}-\eta|$ can be arbitrarily large.

\begin{figure}[h]
	\includegraphics[width=\columnwidth]{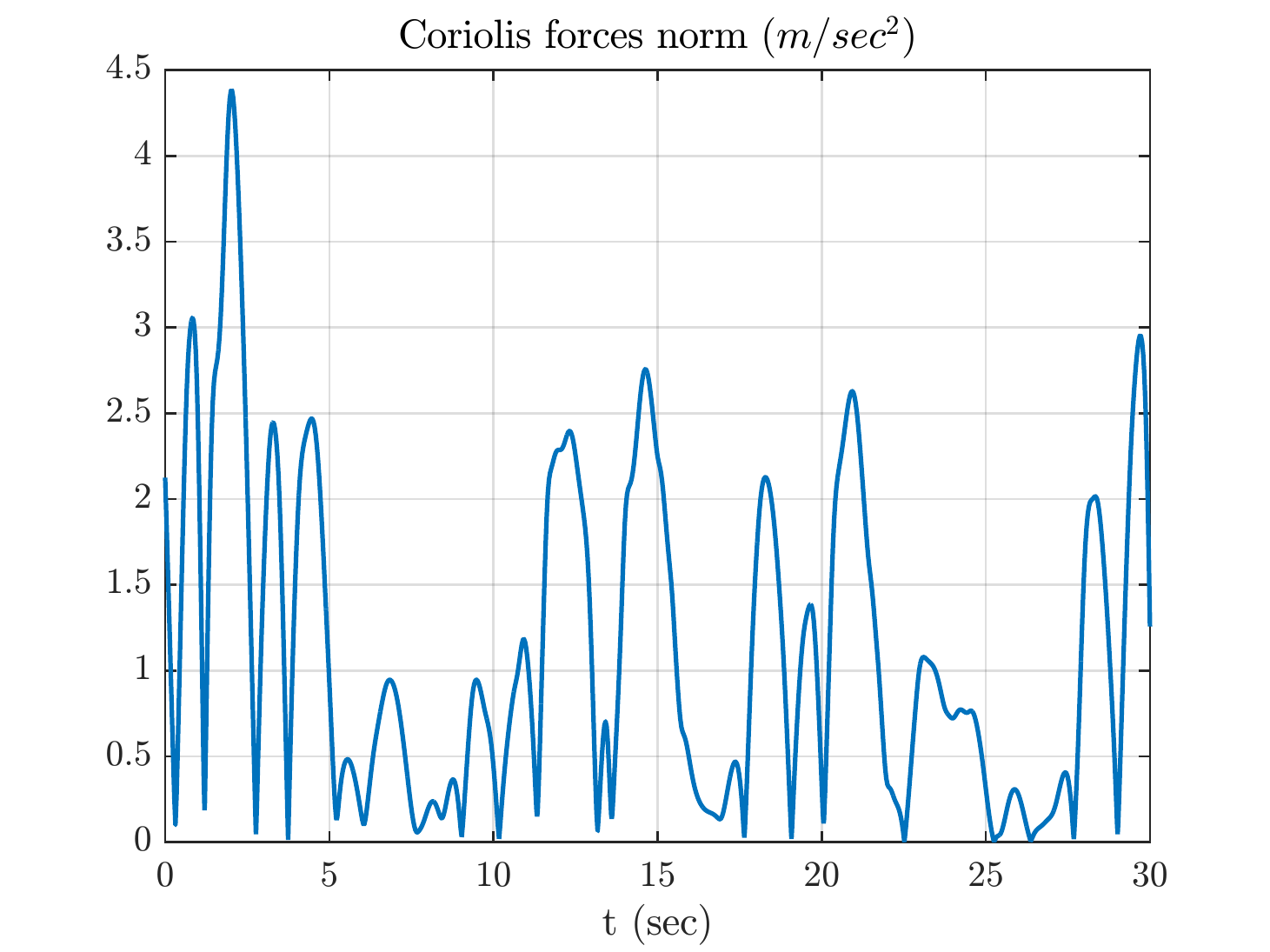}
	\caption{Norm of the Coriolis terms.}
	\label{fig:normC}
\end{figure}

\begin{figure}[h]
	\includegraphics[width=\columnwidth]{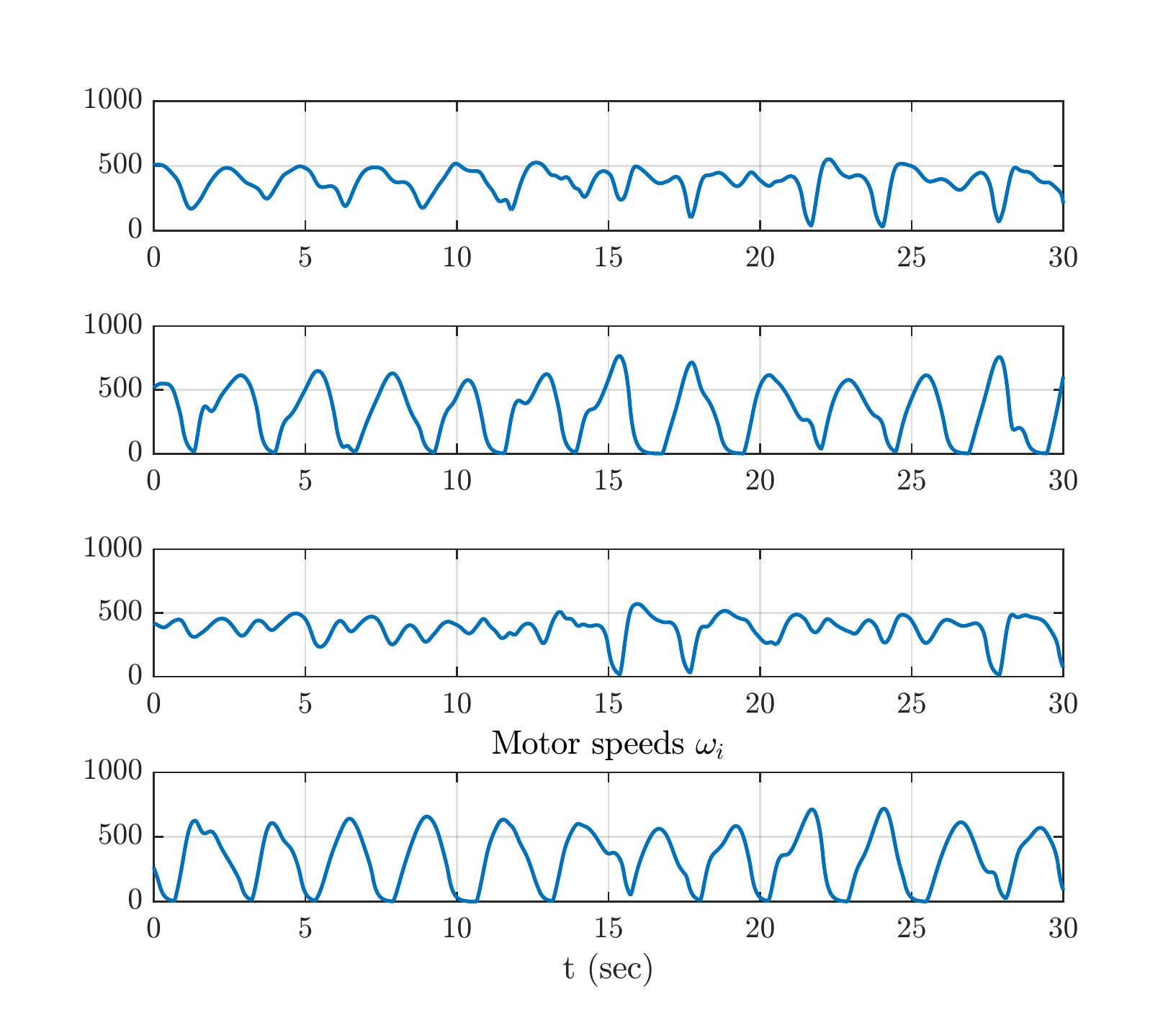}
	\caption{Motor speeds $\omega_1,\omega_2,\omega_3,\omega_4$ ($rad/sec$).}
	\label{fig:Wi}
\end{figure}

\begin{figure}[h]
	\includegraphics[width=\columnwidth]{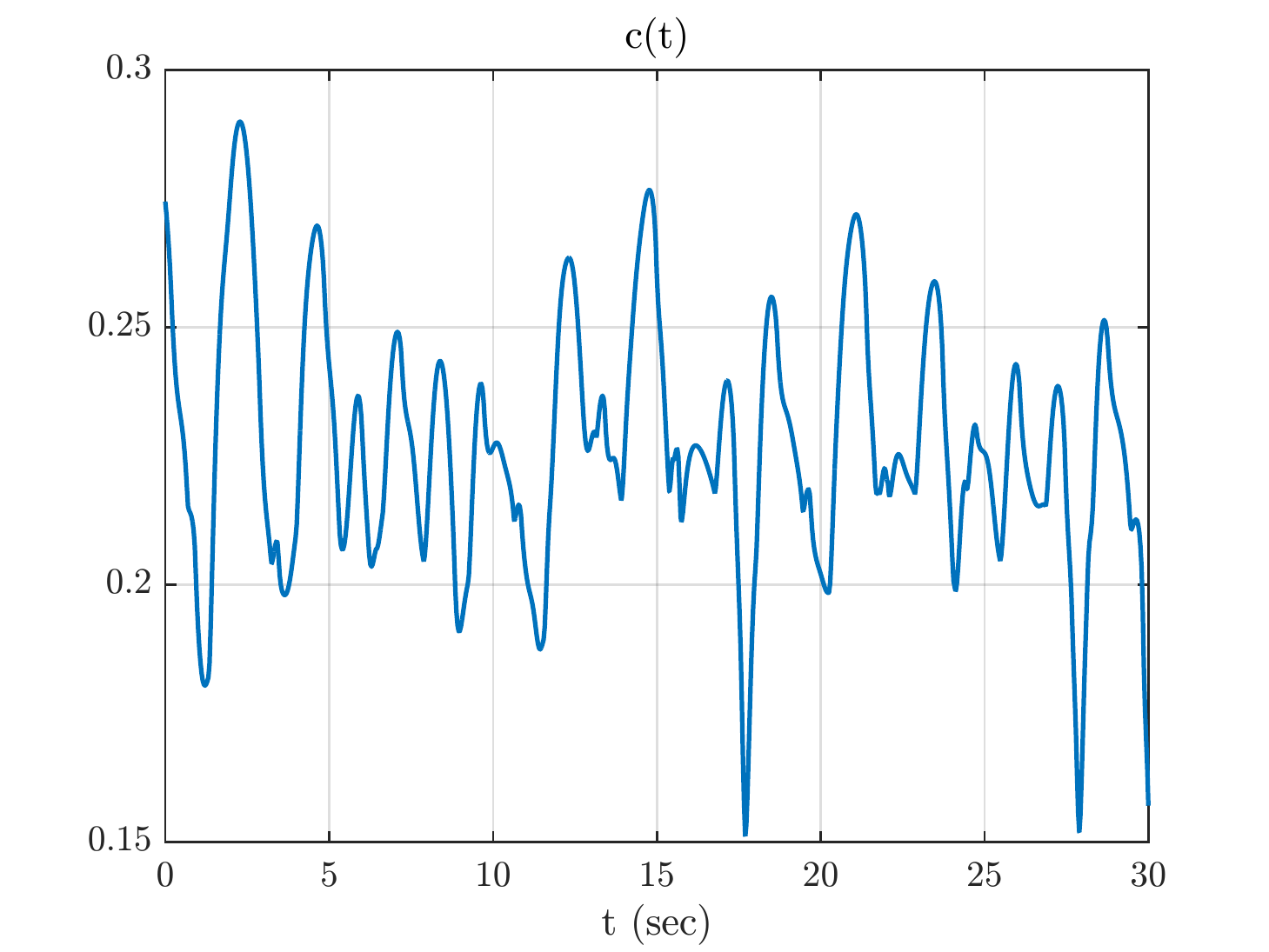}
	\caption{Time evolution of the drag-force coefficient $c(t)$.}
	\label{fig:c}
\end{figure}

We insist that the estimation errors are small, despite the neglected Coriolis forces and imperfect measurements; 
in particular the Coriolis forces are not small, see Fig.~\ref{fig:normC}, since the trajectory is rather aggressive, see Fig.~\ref{fig:Wi}. Finally, nothing significant is lost by using in the observer the nominal value $\bar{c}$ instead of the true value $c(t)$; indeed, the variations of $c$ are quite small, see Fig.~\ref{fig:c}, though the variations of the $\omega_i$ are large as depicted in Fig.~\ref{fig:Wi}. In fact, the main deterioration in the estimates stems from the neglected Coriolis forces rather than from approximating $c$ by $\bar{c}$ or neglecting reasonable measurement imperfections.

Of course, in the case of perfect measurements and with negligible Coriolis terms the observer performs excellently as expected from the theoretical developments. 

\section{Conclusion}\label{sec:conclusions}
We have presented a nonlinear observer for the estimation of the orientation and in-plane velocity of the quadrotor using only the measurements from accelerometers and rate gyros. The design is based on an enhanced model of the quadrotor that includes the drag rotor, and a parametrization that leads to a linear time-varying design model with a nonlinear constraint. The observer has a large (semi-global) domain of convergence and is easily tuned. The robustness of the design with respect to unaccounted modeling (Coriolis) terms, measurement biases and noise, is illustrated in simulation on a rather aggressive trajectory.

%\addtolength{\textheight}{-1cm}   % This command serves to balance the column lengths
                                  % on the last page of the document manually. It shortens
                                  % the textheight of the last page by a suitable amount.
                                  % This command does not take effect until the next page
                                  % so it should come on the page before the last. Make
                                  % sure that you do not shorten the textheight too much.

%________________________________________________________

\bibliographystyle{phmIEEEtran}
\bibliography{bibCDC2016}
%\bibliography{toto}

\end{document}